\documentclass[11pt,a4paper]{article}
\usepackage{amsfonts}
\usepackage[utf8]{inputenc}
\usepackage[english]{babel}
\usepackage{paralist}
\usepackage{url}
\usepackage{dsfont}
\usepackage[dvipsnames]{xcolor}
\usepackage{pgfgantt}
\usepackage{pdflscape}

\usepackage{graphicx}
\usepackage[mathcal]{euscript}
\usepackage{float}
\usepackage{amsmath, amsthm, amssymb}
\usepackage{tikz}
\usepackage{todonotes}
\usepackage{booktabs}
\usepackage{comment}
\usepackage{graphicx}
\usepackage{caption}
\usetikzlibrary{arrows,calc,intersections}

\restylefloat{figure}

\newtheorem{theorem}{Theorem}[section]
\newtheorem{lemma}[theorem]{Lemma}

\theoremstyle{definition}

\theoremstyle{remark}

\usepackage{abstract}
 
\title{Eigenvalue regions and realising monotone stochastic matrices}

\date{}

\author{Brando Vagenende$^{1,2}$, Brecht Verbeken$^{1}$, Marie-Anne Guerry$^{1}$  \\  
        \small $^{1}$Department Business Technology and Operations, Data Analytics Laboratory, Vrije \\
        \small Universiteit Brussel (VUB), Pleinlaan 2, Brussels, 1050, Belgium \\
        \small $^{2}$ Corresponding author: email address: brando.vagenende@vub.be
}

\newlength{\bibitemsep}
\newlength{\bibparskip}
\let\oldthebibliography\thebibliography
\renewcommand\thebibliography[1]{%
  \oldthebibliography{#1}%
  \setlength{\parskip}{\bibitemsep}%
  \setlength{\itemsep}{\bibparskip}%
}
\providecommand{\keywords}[1]
{
  \small	
  \textbf{\textit{Keywords---}} #1
}

\begin{document}
\maketitle

\begin{abstract}
Eigenvalues of stochastic matrices have been studied from two complementary perspectives. The individual eigenvalues are characterised through the well-established Karpelevich regions. The spectrum as a whole has also been analysed, yielding powerful results such as the Johnson–Loewy–London (JLL) inequalities. Current research now turns toward particular subsets of stochastic matrices, among others the doubly stochastic matrices.

This paper studies spectral properties of monotone stochastic matrices which are characterised by the fact that each row stochastically dominates the preceding one, and which arise in contexts such as intergenerational mobility, equal-input models, and credit-rating systems. This paper analyses the dominance matrix associated with a monotone matrix, which is a non-negative matrix that preserves the non-trivial eigenvalues. Properties are established and the conditions are given under which a non-negative matrix can be regarded as a dominance matrix.

In analogy with the stochastic matrices, this study examines for the monotone stochastic matrices both the individual eigenvalues as the spectrum as a whole. Individually, the eigenvalue region for all \(n \times n\) monotone matrices with \(1 \leq n \leq 3\) is completely determined, and realising matrices are provided. Collectively, the set of possible pairs of non-trivial eigenvalues arising from \(3 \times 3\) monotone matrices is characterised, accompanied by realising matrices. In both perspectives, the resulting regions are substantially smaller than those for general stochastic matrices. Finally, this paper proves a reduction theorem stating that, for \(n \geq 4\), the eigenvalue region of \(n \times n\) monotone matrices is contained within that of \((n-1) \times (n-1)\) stochastic matrices.
\end{abstract}

\keywords{non-negative matrices, stochastic matrices, monotone stochastic matrices, eigenvalues, spectrum, eigenvalue regions}

\section{Introduction}

A stochastic matrix is a non-negative matrix with each row sum equal to one. These matrices find their way into numerous applications in different domains, most famously in Markov chain theory \cite{Bartholomew}. In this field, stochastic matrices describe the transition probabilities of a given process over time. The spectrum of a stochastic matrix provides insights for the aforementioned Markov chains, among others the steady-state behaviour, the long-term behaviour and the convergence 
properties \cite{pillai2005perron, delbianco2023markov, racoceanu1995new, meyer2000applied}. Characterising the eigenvalue regions of stochastic matrices is a challenging problem with a long and extensive research history. The eigenvalues of such matrices can be examined from distinct perspectives.

On the one hand, the set of the individual eigenvalues can be considered. It was as early as 1938 that Kolmogorov first came up with this problem. Specifically, he examined the set $\mathcal{S}_n$ of $n \times n$ stochastic matrices and asked whether an exact description could be found for the region \(\Theta_n = \bigl\{ \lambda \in \mathbb{C} : \lambda \text{ is an eigenvalue of an } n \times n \text{ stochastic matrix} \bigr\}\), where $n$ is a fixed natural number. It was only in 1946 that Dmitriev and Dynkin \cite{Dynkin} could give part of the solution to this problem, namely a complete description of the eigenvalue region for $\mathcal{S}_n$ up to $n=5$. In 1951, Karpelevich \cite{Karpelevich} was finally able, by generalising some ideas of Dmitriev and Dynkin, to give a complete description of the region consisting of all eigenvalues for $\mathcal{S}_n$ for any $n$.

On the other hand, eigenvalues can be considered collectively. In 1978, Loewy and London \cite{loewy1978note} investigated when a specified set of complex numbers can be realised as the spectrum of a non-negative matrix. They introduced the well-known JLL-inequalities as necessary conditions on the spectrum, which are formulated in terms of traces of matrix powers. For the specific case of \(3 \times 3\) matrices, they presented additional conditions under which such spectra are realisable. The non-negative inverse eigenvalue problem remains an active area of research, with recent developments including those in e.g. \cite{johnson2025perron}.

Within this framework, doubly stochastic matrices, a subset of stochastic matrices with both row and column sums equal to one, are also a common topic of discussion. Research on the eigenvalue regions of these matrices is ongoing, with several results established, though much about the eigenvalues remains unknown. Conjectures and partial proofs appear in several publications, e.g. \cite{mashreghi2007conjecture}. Besides (doubly) stochastic matrices, the eigenvalue regions of other matrix types, such as Metzler matrices \cite{domka2022spectrum} and Leslie matrices \cite{kirkland1992eigenvalue}, are also studied.

This research focuses specifically on a subset of the \(n \times n\) stochastic matrices \(\mathcal{S}_n\), namely the set $\mathcal{M}_n$ of $n \times n$ monotone matrices. This terminology is not always used consistently in previous work. We adopt the definition provided by Daley \cite{daley1968stochastically}: A monotone matrix \(M =(m_{ij})\) is a stochastic matrix in which each row is stochastically dominated by the next row, i.e. \(\sum_{j=r}^n m_{lj} \geqslant \sum_{j=r}^n m_{kj} \forall l > k, \hspace{2mm} \forall r \in \{1, \ldots, n \} .\) Monotone matrices appear in various contexts, such as intergenerational occupational mobility \cite{conlisk1990monotone}, equal-input modeling \cite{baake2022equal} and credit ratings based 
systems \cite{jarrow1997markov}. While some research has been conducted on these matrices and their eigenvalues, such as in \cite{guerry2022monotone}, the scope remains relatively narrow. 

On the one hand, by examining eigenvalues individually, this work provides new insights into the eigenvalue region
\[
\Xi_{n} := \bigl\{\,\lambda \in \mathbb{C} \,\bigm|\, \exists\, M \in \mathcal{M}_{n} \text{ such that } \lambda \in \sigma(M)\,\bigr\},
\]
defined for arbitrary $n$, where $\sigma(M)$ denotes the spectrum of $M$.
Similar to the (doubly) stochastic matrices, it also holds for the monotone matrices that \(\Xi_n \subset \Xi_{n+1}\) for every \(n \geq 1\). On the other hand, by examining the eigenvalues collectively, and knowing that the spectrum of a $3 \times 3$ monotone matrix is real \cite{guerry2022monotone}, the set
\[
\xi_{3} := \bigl\{\,(\lambda_{2}, \lambda_{3}) \in \mathbb{R}^{2} \,\bigm|\,
\exists\, M \in \mathcal{M}_{3} \text{ such that } \sigma(M)=\{1,\lambda_{2},\lambda_{3}\}
\text{ with } \lambda_{2} \ge \lambda_{3} \,\bigr\}
\]
is fully characterised for $n=3$. 

This paper presents, in Section \ref{section_dominance_matrix}, several properties of the dominance matrix $D(M)$ associated with a monotone matrix $M$, together with an explicit construction. Section \ref{Monotone eigenvalue regions} provides a full determination of $\Xi_{1}$, $\Xi_{2}$, and $\Xi_{3}$, along with corresponding realising matrices. In Section \ref{eigenvalue_pairs}, the set $\xi_{3}$ is fully characterised, and accompanied by realising matrices. Section \ref{Reduction theorem} proves a reduction theorem which gives insight into the eigenvalue region for $\mathcal{M}_n$ for $n \geqslant 4$. Finally, in Section \ref{Conclusions and further research}, conclusions and further research avenues are presented.

\section{Dominance matrix}\label{section_dominance_matrix}

To analyse the eigenvalues of an \(n \times n\) monotone matrix \(M\), the dominance matrix \(D(M)\) plays a central role. Its key advantage is that it reduces the problem to a, one order lower, non-negative $(n-1)\times(n-1)$ matrix while retaining exactly the same nontrivial eigenvalues, i.e \(\sigma(D(M))=\sigma(M)\setminus\{1\}\). The dominance matrix \(D(M)\) is given by \[(D(M))_{kl} = \sum_{j=1}^l m_{kj} - \sum_{j=1}^l m_{k+1,j}, \forall k,l \in \{1, \ldots, n-1 \}.\] By the definition of a monotone matrix, it follows directly that $D(M)\geq 0$. More details can be found in \cite{conlisk1990monotone}.

It is natural to ask whether the dominance matrix can be an arbitrary $(n-1)\times(n-1)$ non-negative matrix. However, the additional constraints inherent to monotone matrices imply that this is not the case. Lemma \ref{properties_dominance_matrix} proves necessary conditions for the case $n=3$.

Given this, let $M=(m_{ij})$ be a $3\times3$ monotone matrix. Its associated dominance matrix is
\[
D(M)=
\begin{pmatrix}
    a & b \\
    c & d
\end{pmatrix},
\]
where
\begin{align}
    a &= m_{11} - m_{21}, \label{a}\\
    b &= (m_{11} + m_{12}) - (m_{21} + m_{22}), \label{b}\\
    c &= m_{21} - m_{31}, \label{c}\\
    d &= (m_{21} + m_{22}) - (m_{31} + m_{32}). \label{d}
\end{align}
The nontrivial eigenvalues of $M$ are precisely the eigenvalues of $D(M)$, given explicitly by
\[
\lambda_2 =
\frac{a + d + \sqrt{(a-d)^2 + 4bc}}{2}
\quad\text{and}\quad
\lambda_3 =
\frac{a + d - \sqrt{(a-d)^2 + 4bc}}{2}.
\]

Further, the following properties can be deduced.

\begin{lemma}\label{properties_dominance_matrix}
Let $D(M)$ be the dominance matrix defined above. Then:
\begin{enumerate}
    \item Column- and anti-diagonal-sums are bounded by 1:
    \[
        a + c \le 1, \qquad
        b + c \le 1, \qquad
        b + d \le 1.
    \]

    \item Column- and anti-diagonal-products are bounded by \(\frac{1}{4}\):
    \[
        ac \le \tfrac{1}{4}, \qquad
        bc \le \tfrac{1}{4}, \qquad
        bd \le \tfrac{1}{4}.
    \]

    \item Non-negative trace:
    \[
        \mathrm{tr}[D(M)] = a + d \ge 0.
    \]

    \item Determinant bound:
    \[
        \det[D(M)] = ad - bc \ge -\tfrac{1}{4}.
    \]
\end{enumerate}
\end{lemma}
\begin{proof}
\begin{enumerate}
    \item Using (\ref{b}) and (\ref{c}), we obtain
    \(b + c = (m_{11} + m_{12}) - (m_{22} + m_{31}) \leq m_{11} + m_{12} \leq 1,
    \) which establishes \(b+c \leq 1\). The arguments for \(a+c \leq 1\) and \(b+d \leq 1\) proceed in the same manner.
    \item According to (\ref{b}) an upper bound for \(b\) can be obtained by choosing \(m_{22} = 0\) and \(m_{11} + m_{12} = 1\). In turn, according to (\ref{c}), \(c\) can be bounded upwards by setting \(m_{31} = 0\). In this way, we get \( b \cdot c \leq (1-m_{21}) \cdot m_{21}.\) The maximum of this last expression is reached for \(m_{21}=1/2\) and is equal to \(1/4\). So we obtain \(bc \leq \frac{1}{4}\). Similar arguments provide \(ac \leq \frac{1}{4}\) and \(bd \leq \frac{1}{4}\).
    
    \item This follows from the fact that \(D(M) \geq 0\) (see \cite{guerry2022monotone}).
    \item From \(bc \le \tfrac{1}{4}\) (proved in 2) follows that \(\mathrm{det}[D(M)]=ad - bc \ge -bc \ge -\tfrac{1}{4}\).
\end{enumerate}
\end{proof}

The properties of Lemma \ref{properties_dominance_matrix} demonstrate that not every random non-negative $2\times 2$ matrix
\(\begin{pmatrix} a & b \\ c & d\end{pmatrix}
\)
is the dominance matrix of a $3\times 3$ monotone matrix. However, the  inequalities in Lemma \ref{conditions dominance and monotone} form a necessary and sufficient condition:
\begin{lemma}\label{conditions dominance and monotone}
A non-negative matrix $\begin{pmatrix} a & b \\ c & d \end{pmatrix}$ arises as the dominance matrix \newline $D(M)$ of a $3\times 3$ monotone matrix $M$ if and only if there exist $m_{11}, m_{33} \in [0,1]$ such that
\[
\begin{cases}
    a + c \le m_{11}, \\
    b + d \le m_{33}, \\
    m_{11} + m_{33} \le 1 + b + d, \\[4pt]
    m_{11} + m_{33} \le 1 + a + d, \\[4pt]
    m_{11} + m_{33} \le 1 + a + c.
\end{cases}
\]
\end{lemma}

\begin{proof}
\(\Longrightarrow\)
Assume that \(M=(m_{ij})\) is a \(3\times3\) monotone matrix with
\(D(M)=
\begin{pmatrix}
    a & b \\
    c & d
\end{pmatrix}.
\)
Using \eqref{a}--\eqref{d}, together with the stochasticity of \(M\), i.e.,
\(m_{11}+m_{12}=1-m_{13}\), \(m_{21}+m_{22}=1-m_{23}\), and
\(m_{31}+m_{32}=1-m_{33}\), we obtain
\begin{align}
    m_{21} &= m_{11} - a, \label{m21} \\
    m_{31} &= m_{11} - a - c, \label{m31} \\
    m_{13} &= m_{33} - b - d, \label{m13} \\
    m_{23} &= m_{33} - d. \label{m23}
\end{align}
Since \(M \ge 0\), relations \eqref{m31} and \eqref{m13} immediately yield
the first two inequalities in the system. Moreover, from
\(m_{12}=1-m_{11}-m_{13}\) and \eqref{m13}, we deduce
\[
m_{11} + m_{33} \le 1 + b + d,
\]
which provides one of the remaining inequalities. The remaining two follow
by analogous considerations.

\(\Longleftarrow\)
Conversely, suppose that the inequalities in the statement of the lemma are satisfied. Then one may explicitly construct the following
$3\times3$ monotone matrix:
\[
M =
\begin{pmatrix}
    m_{11} 
    & 1 - (m_{11} + m_{33} - b - d)
    & m_{33} - b - d \\[6pt]
    m_{11} - a
    & 1 - (m_{11} - a + m_{33} - d)
    & m_{33} - d \\[6pt]
    m_{11} - a - c
    & 1 - (m_{11} - a - c + m_{33})
    & m_{33}
\end{pmatrix}
\]
with dominance matrix
\(
D(M)=\begin{pmatrix} a & b \\ c & d\end{pmatrix}.
\)
\end{proof}

Observe that each monotone matrix determines a unique dominance matrix, but different monotone matrices can have the same dominance matrix. For example, for the monotone matrices \(M_1=\begin{pmatrix}
    0,3 & 0,7 & 0 \\
    0,2 & 0,7 & 0,1 \\
    0,1 & 0,7 & 0,2
\end{pmatrix}\) and \(M_2=\begin{pmatrix}
    0,4 & 0,6 & 0 \\
    0,3 & 0,6 & 0,1 \\
    0,2 & 0,6 & 0,2
\end{pmatrix}\) holds that \(D(M)=\begin{pmatrix}
    0,1 & 0,1 \\
    0,1 & 0,1
\end{pmatrix}\). This is not an issue for our purposes, since we are concerned only with the existence of a monotone matrix possessing a prescribed eigenvalue, rather than with the multiplicity of such matrices.

\section{Eigenvalue regions $\Xi_1, \Xi_2$ and \(\Xi_3\)}\label{Monotone eigenvalue regions}

In analogy with the Karpelevich regions \(\Theta_n\) \cite{Karpelevich}, we consider the eigenvalues individually in this section and determine the eigenvalue regions \(\Xi_n\) for \(1 \leqslant n \leqslant 3\). Because monotonicity imposes further restrictions on the matrix entries, we expect the resulting eigenvalue regions to be narrower than those that arise for the stochastic matrices.

\begin{theorem}
    \( \Xi_1 = \{ 1 \} \).
\end{theorem}
\begin{proof}
    This statement is trivial since a \(1 \times 1\) monotone matrix can only be the matrix \((1)\) with eigenvalue \(\lambda=1\).
\end{proof}
In this trivial case, it holds that \(\Xi_1 = \Theta_1\).
\begin{theorem}
    \( \Xi_2 = [0,1] \).
\end{theorem}
\begin{proof}
Let \(M\) be a \(2 \times 2\) monotone matrix with \(\sigma(M)= \{\lambda_1 = 1, \lambda_2 \}\), then \(\lambda_2 \geqslant 0\) because of \(\mathrm{tr}(M) \geqslant 1\) (see \cite{guerry2022monotone}). Hence, \(\Xi_2 \subseteq [0,1]\). Furthermore, for \(\lambda_2 \in [0,1] \) the monotone matrix \(\begin{pmatrix}
                \lambda_2 & 1 - \lambda_2 \\
                0 & 1 
            \end{pmatrix}\) has \(\lambda_2\) as eigenvalue. This concludes the proof.
\end{proof}
We get, in this case, a first reduction, namely that \(\Xi_2=[0,1]\subset [-1, 1]=\Theta_2\).
\vspace{15 mm}
\begin{theorem}\label{Xi_3}
    \(\Xi_3 = [-1/2, 1]\).
\end{theorem}
\begin{proof}

\textbf{Step 1: proof of \(\Xi_3 \subseteq [-1/2, 1]\)}.
\vspace{3mm}

Let \(M\) be an \(3 \times 3\) monotone matrix and \(D(M)\) its dominance matrix. In order to prove \(\Xi_3 \subseteq [-1/2, 1]\), we want to know how small an eigenvalue can be. Therefore, we are going to minimise the smallest eigenvalue (see Section \ref{section_dominance_matrix}):
    \begin{align}
        \lambda_3 &= \frac{a + d - \sqrt{(a-d)^2 + 4bc}}{2} \nonumber \\
         &\geqslant \frac{a + d - \sqrt{(a+d)^2 + 4bc}}{2} \label{ineq1}\\
         &\geqslant \frac{a + d - \sqrt{(a+d + 2 \sqrt{bc})^2}}{2} \label{ineq2} \\
         &= - \sqrt{bc} \label{eq1}
    \end{align}
    
    Above, inequalities (\ref{ineq1}) and (\ref{ineq2}) follow from the fact that \(a\) and \(d\) are positive. It follows from equation (\ref{eq1}) that, in order to investigate a lower bound for the eigenvalues, we need to maximise the expression \(b \cdot c\). From Lemma \ref{properties_dominance_matrix}(2) follows that this maximum is equal to \(1/4\). So we obtain that \( \lambda_3 \geqslant - 1/2 \), and as we know that an eigenvalue of a stochastic matrix is at most 1, this concludes the first step.

\vspace{3mm}
\textbf{Step 2: \(\Xi_3 = [-1/2, 1]\)}
\vspace{3mm}

So we already know that the region \(\Xi_3\) is located in the interval \([-1/2, 1]\). However, the last interval turns out to be exactly the monotone eigenvalue region \(\Xi_3 \). This can be seen by the following 2 constructions of realising matrices.

\vspace{3mm}
\textbf{Realising matrices type 1 - covering \([0,1]\)}
\vspace{3mm}

For each \( \alpha \in [0,1]\), we construct the following monotone matrix

\[\begin{pmatrix}
    0 & 1-\alpha & \alpha \\
    0 & 1-\alpha & \alpha \\
    0 & 0 & 1
\end{pmatrix} .\]

Via a simple calculation, we obtain the following characteristic equation

\[(\lambda - 1)(-\lambda^2 - \lambda \alpha + \lambda) = 0,\]

from which the eigenvalues below follow:
\[ \lambda_1 = 1, \lambda_2 = 1-\alpha \text{ and } \lambda_3 = 0  .\]

For \( \alpha \in [0,1] \), the eigenvalue \( \lambda_2 = 1 - \alpha \) traverses the line segment \([0,1]\).

\vspace{3mm}
\textbf{Realising matrices type 2 - covering \([-1/2,0]\)}
\vspace{3mm}

For each \( \alpha \in [0, 1/2]\), we construct the following monotone matrix

\[ \begin{pmatrix}
    1/2 - \alpha & 1/2 + \alpha & 0 \\
    1/2 - \alpha & 0 & 1/2 + \alpha \\
    0 & 1/2 - \alpha & 1/2 +\alpha
\end{pmatrix}\]

with characteristic equation

\[-\lambda^3 + \lambda^2 + (1/4 - \alpha^2) \lambda + (\alpha^2 - 1/4) = 0,\]

from which the eigenvalues below follow:
\[ \lambda_1 = 1,  \lambda_2 = \sqrt{1/4 - \alpha^2} \text{ and } \lambda_3 = - \sqrt{1/4 - \alpha^2}.\]

For \( \alpha \in [0,\frac{1}{2}] \), the eigenvalues \(\lambda_2 = \sqrt{1/4 - \alpha^2}\) and \(\lambda_3 = - \sqrt{1/4 - \alpha^2} \) traverse respectively line segments \([-1/2,0]\) and \([0, 1/2]\).
\end{proof}

\begin{figure}[ht]
    \centering
\includegraphics[width=0.4\textwidth]{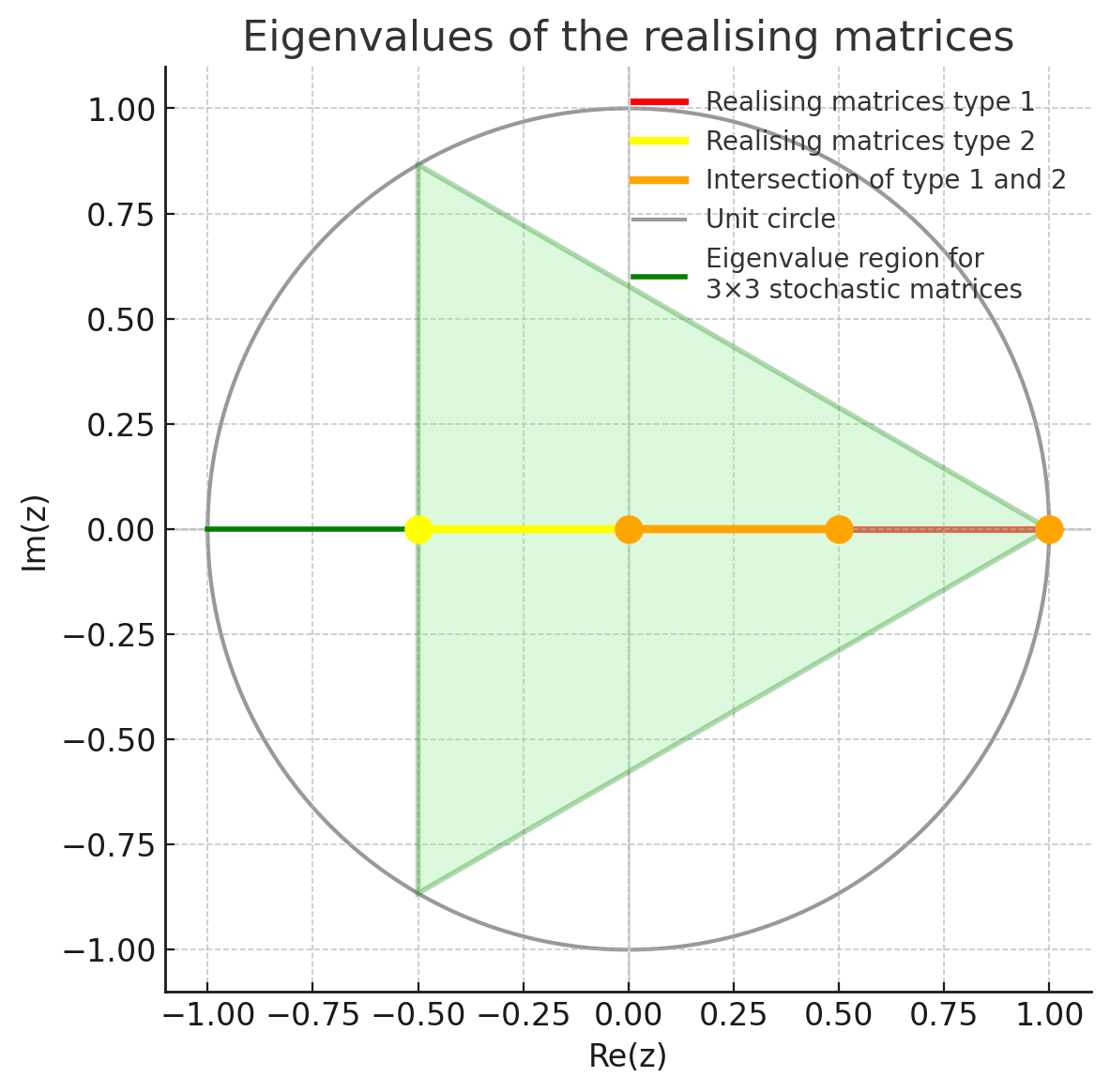}
    \caption{Eigenvalues of the realising matrices.}
    \label{eigenvalues of realising matrices}
\end{figure}

The eigenvalues of the realising matrices of type 1 and 2 are presented in Figure \ref{eigenvalues of realising matrices} where the union results in the eigenvalue region \(\Xi_3 =[-1/2, 1]\). Figure \ref{eigenvalues of realising matrices} also illustrates the eigenvalue region $\Theta_3$ for $3 \times 3$ stochastic matrices, namely
\(
\Theta_3 = [-1, \tfrac{1}{2}) \,\cup\, \mathrm{conv}\{1, e^{\frac{i2\pi}{3}}, e^{\frac{i4\pi}{3}}\},
\)
where $\mathrm{conv}\{\cdot\}$ denotes the convex hull. We see that \(\Xi_3\) represents a significant contraction of \(\Theta_3\). For instance, in contrast to \(\Theta_3\), \(\Xi_3\) does not consist of complex eigenvalues. It is worth noting that this is a remarkable feature, as the eigenvalue region for \(3 \times 3\) doubly stochastic matrices coincides with that of stochastic matrices and therefore does not exhibit such a reduction.

\section{Possible eigenvalue pairs \((\lambda_2, \lambda_3)\) for \(\mathcal{M}_3\)}\label{eigenvalue_pairs}

In analogy with the JLL-inequalities \cite{loewy1978note}, we investigate in this section the set \(\xi_3\) of possible eigenvalue pairs $(\lambda_{2}, \lambda_{3})$ arising from a $3\times3$ monotone matrix \(M\). 
Since $\lambda_{1}=1$ holds universally, it is excluded from consideration. 
Although the previous section established that all eigenvalues of a $3\times3$ monotone matrix must lie in the interval $[-1/2,\,1]$, this bound does not identify which pairs can occur simultaneously. For instance, $\lambda_2 = 1$ and $\lambda_3 = -\frac{1}{2}$ both belong to $[-\frac{1}{2}, 1]$ but their combination is not possible, since $\lambda_2 \lambda_3 = -\tfrac{1}{2}$, which contradicts Lemma \ref{properties_dominance_matrix}(4). Our aim is therefore to characterise the region \(\xi_3\) in the $(\lambda_{2}, \lambda_{3})$-plane. By similar arguments to \cite{vagenende2025star}, we get:
\begin{lemma}\label{starconvexity_of_xi_3}
    The region \(\xi_3\) is star-convex with respect to the origin. 
\end{lemma}
Thus, determining the boundary of \(\xi_3\) is sufficient to obtain the complete region.

Let
\(
D(M) = 
\begin{pmatrix}
    a & b \\
    c & d
\end{pmatrix}
\)
be the dominance matrix as introduced in Section \ref{section_dominance_matrix}. For particular subsets of $\lambda_{2}$ and $\lambda_{3}$, additional information about $\mathrm{tr}[D(M)]$ can be derived.
\begin{lemma}\label{hulplemma}
If \(\lambda_2 \geq \frac{1+\sqrt{5}}{4} \) and \(\lambda_3 \leq 0\), then \(\frac{1}{2} \leq \mathrm{tr}[D(M)]=a+d \leq 1 \).
\end{lemma}
\begin{proof}
Because $\lambda_{2} \le 1$ and $\lambda_{3} \le 0$, we immediately obtain $\mathrm{tr}[D(M)] = \lambda_{2} + \lambda_{3} \le 1$. Moreover, for 
$\lambda_{2} \ge \frac{1+\sqrt{5}}{4}$, the smallest admissible value of $\lambda_{3}$ is $-\tfrac{1}{\,1+\sqrt{5}\,}$, since Lemma \ref{properties_dominance_matrix}(4) implies $\mathrm{det}[D(M)]=\lambda_{2}\lambda_{3} \ge -\frac{1}{4}$. Evaluating the trace at these extremal values yields \(\frac{1+\sqrt{5}}{4} - \frac{1}{1+\sqrt{5}} = \frac{1}{2}\), and therefore $\mathrm{tr}[D(M)] \ge \frac{1}{2}$ throughout this region.
\end{proof}
In fact, a more explicit connection can be derived between $\mathrm{tr}[D(M)]$ and $\det[D(M)]$.
\begin{lemma}\label{Link tr det}
    If \(\lambda_2 \geq \frac{1+\sqrt{5}}{4}\) and \(\lambda_3 \leq 0\), then \[\mathrm{tr}[D(M)]^2 - \mathrm{tr}[D(M)] - \mathrm{det}[D(M)] \leq 0.\]
\end{lemma}
\begin{proof}
 Substituting $\mathrm{tr}[D(M)] = a + d$ and $\det[D(M)] = ad - bc$ into $\mathrm{tr}[D(M)]^{2} - \mathrm{tr}[D(M)] - \det[D(M)] \le 0$ and rearranging terms results in the equivalent inequality 
\begin{equation}\label{inequality tr and det}
ad + bc - a(1 - a) - d(1 - d) \le 0.
\end{equation}
Using $d \le 1 - a$ and $b \le 1 - d$, as established in Lemma \ref{properties_dominance_matrix}(1), we conclude that inequality \eqref{inequality tr and det} holds whenever $c \le d$. Consequently, the remainder of the proof may focus on the case $c > d$.

By the bound $b \le 1 - c$ (from Lemma \ref{properties_dominance_matrix}(1)), we obtain
\[ad + bc - a(1 - a) - d(1 - d) \le G(c)\]
where
\(G(c) := a^{2} + d^{2} + ad - a - d + c - c^{2}\) is considered for fixed $a$ and $d$. Since we are in the case $c > d$ and Lemma \ref{properties_dominance_matrix}(1) ensures $c \le 1 - a$, the admissible range is $c \in ]d,\,1 - a]$. Our goal is therefore to show that $G(c) \le 0$ for all 
$c \in ]d,\,1 - a]$.

 Analysing $G(c)$ reduces to studying the function $c \mapsto c-c^2$, which is concave with a unique maximum at $c = \tfrac{1}{2}$, increasing on $[0,\,\tfrac{1}{2})$ and decreasing on $(\tfrac{1}{2},\,1]$. We distinguish three subcases.

In case of $1 - a \le \tfrac{1}{2}$, it holds that
\[G(c) \leq G(1-a) = d(a + d - 1) \leq 0,\]
since $d \ge 0$ and Lemma \ref{hulplemma} ensures $a + d \le 1$.

In case of $1 - a > \tfrac{1}{2}$ and $d \ge \tfrac{1}{2}$, we obtain \[G(c) \leq G(d) = a(a + d - 1) \leq 0.\]

In case of \(d < \frac{1}{2} < 1-a\), it can be observed that
\[G(c) \leq G\left(\frac{1}{2}\right)= a^2 -a + d^2 - d + ad + \frac{1}{4}:= P(a,d).\]
Since \(1-a > \frac{1}{2}\), it follows that \(a < \frac{1}{2}\). Combined with \(d < \frac{1}{2}\) and Lemma \ref{hulplemma}, this implies that the domain of \(P(a,d)\) is the triangular region bounded by \(a < \frac{1}{2}\), \(d < \frac{1}{2}\), and \(a + d \ge \frac{1}{2}\). The function \(P(a,d)\) possesses a single extremum, which is a minimum at \(\bigl(\frac{1}{3}, \frac{1}{3}\bigr)\). It therefore follows that the function is bounded above by its boundary values
\[P\left(a, \frac{1}{2}\right) = P\left(a, \frac{1}{2}-a\right)= a\left(a-\frac{1}{2}\right) \leqslant 0 \text{ and } P\left(\frac{1}{2}, d\right)=d\left(d-\frac{1}{2}\right) \leqslant 0,\]
from which follows that \(G(c) \leqslant P(a,d) \leqslant0\).

From the three subcases above, we conclude that $G(c) \le 0$ for all 
$c \in (d,\,1-a]$, which completes the proof.
\end{proof}
The following theorem provides a complete characterization of $\xi_3$ in the $(\lambda_2,\lambda_3)$-plane, as illustrated in Figure \ref{fig:possible1}.

\noindent
\begin{minipage}{0.45\textwidth}
    \centering
    \includegraphics[width=\textwidth]{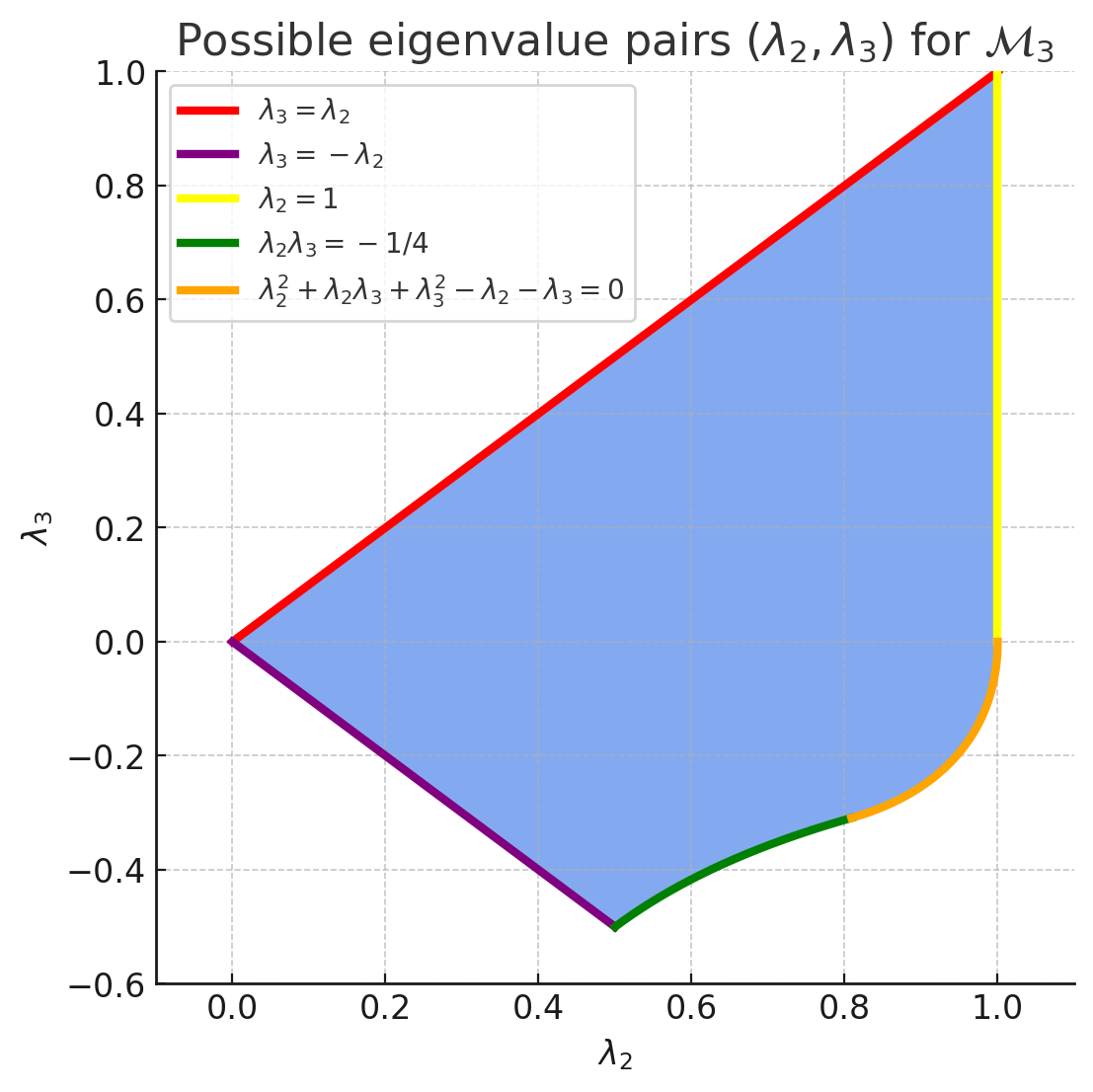}
    \captionof{figure}{Possible eigenvalue pairs \((\lambda_2, \lambda_3)\) for \(\mathcal{M}_3\).}
    \label{fig:possible1}
\end{minipage}
\hfill
\begin{minipage}{0.50\textwidth}
    \centering
    \includegraphics[width=\textwidth]{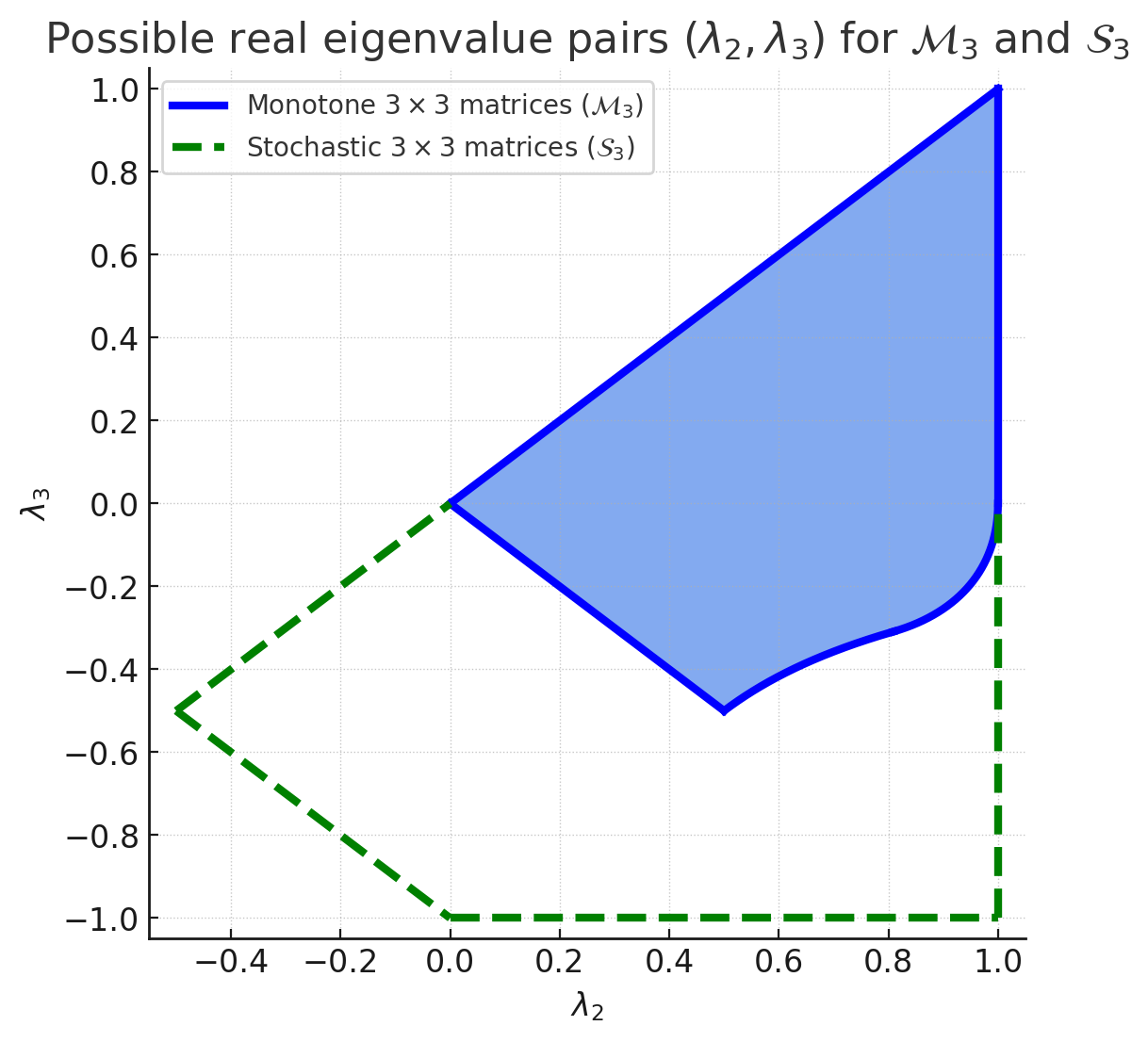}
    \captionof{figure}{Possible real eigenvalue pairs \((\lambda_2, \lambda_3)\) for \(\mathcal{M}_3\) and \(\mathcal{S}_3\).}
    \label{fig:possible2}
\end{minipage}
\vspace{5 mm}
\begin{theorem}
    \(\xi_3\) is the region bounded by the following curves:
    \begin{itemize}
        \item \(C_1:\) \(\lambda_2 = \lambda_3 \) for \(0 \leq \lambda_2 \leq 1\),
\item \(C_2:\) \(\lambda_2 = - \lambda_3 \) for \(0 \leq \lambda_2 \leq \frac{1}{2}\),
\item \(C_3:\) \(\lambda_2 = 1\) for \(0 \leq \lambda_3 \leq 1 \),
\item \(C_4:\) \(\lambda_2\cdot\lambda_3 = -\frac{1}{4} \) for \(\frac{1}{2} \leq \lambda_2 \leq \frac{1+\sqrt{5}}{4}\),
\item \(C_5:\) \(\lambda_2^2 + \lambda_2 \lambda_3 + \lambda_3^2 - \lambda_2 - \lambda_3 \leq 0 \) for \(\frac{1+\sqrt{5}}{4} \leq \lambda_2 \leq 1\) and \(\lambda_3 \leq 0\).
    \end{itemize}
\end{theorem}
\begin{proof}
By Lemma \ref{starconvexity_of_xi_3} it suffices to determine the boundary of \(\xi_3\). We first establish the boundary curves and subsequently show that each of them is realised. Thus, we begin by identifying the boundaries, which are given below:
    \begin{itemize}
        \item $C_1$ follows directly from the imposed condition $\lambda_2 \ge \lambda_3$.
        \item From Lemma \ref{properties_dominance_matrix}(3) we know that \(\mathrm{tr}[D(M)]=\lambda_2 + \lambda_3 \geq 0\), which gives \(C_2\).
        \item Since for a stochastic matrix, the modulus of each of the eigenvalues is at most 1, \(C_3\) follows trivially.
        \item From Lemma \ref{properties_dominance_matrix}(4) follows immediately that \(\mathrm{det}[D(M)]= \lambda_2 \lambda_3 \geq - \frac{1}{4}\), and hence \(C_4\).
        \item For $\lambda_2 \ge \frac{1+\sqrt{5}}{4}$ and $\lambda_3 \le 0$, Lemma \ref{Link tr det} yields $\mathrm{tr}[D(M)]^{2} - \mathrm{tr}[D(M)] - \mathrm{det}[D(M)] \le 0$. Substituting $\mathrm{tr}[D(M)] = \lambda_1 + \lambda_3$ and $\mathrm{det}[D(M)] = \lambda_1 \lambda_3$ into this expression gives $C_5$.
    \end{itemize}
Combining these boundary conditions shows that the region is precisely bounded by the curves \(C_i\) specified in this theorem. This is illustrated in Figure \ref{fig:possible1}.

What remains is to verify that each of these bounds is realised. This can be confirmed by the matrices presented in Table \ref{table}.
\end{proof}

\begin{table}[h]
\centering
\renewcommand{\arraystretch}{1.4}
\begin{tabular}{c c c}
\toprule
\textbf{Curve} & \textbf{Realising Matrices} & \textbf{Parameter Range} \\
\midrule

$C_1$ &
$\begin{pmatrix}
    \frac{1+2\alpha}{3} & \frac{1-\alpha}{3} & \frac{1-\alpha}{3} \\
    \frac{1-\alpha}{3} & \frac{1+2\alpha}{3} & \frac{1-\alpha}{3} \\
    \frac{1-\alpha}{3} & \frac{1-\alpha}{3} & \frac{1+2\alpha}{3}
\end{pmatrix}$ 
& $\alpha \in [0,1]$ \\

$C_2$ &
$\begin{pmatrix}
    \alpha & 1-\alpha & 0 \\
    \alpha & 1-2\alpha & \alpha \\
    0 & 1-\alpha & \alpha
\end{pmatrix}$
& $\alpha \in \left[0,\frac{1}{2}\right]$ \\

$C_3$ &
$\begin{pmatrix}
    \alpha & 1-\alpha & 0 \\
    0 & 1 & 0 \\
    0 & 0 & 1
\end{pmatrix}$
& $\alpha \in [0,1]$ \\

$C_4$ &
$\begin{pmatrix}
    1-\alpha & \alpha & 0 \\
    \frac{1}{2} & 0 & \frac{1}{2} \\
    0 & \frac{1}{2} & \frac{1}{2}
\end{pmatrix}$
& $\alpha \in \left[0,\frac{1}{2}\right]$ \\

$C_5$ &
$\begin{pmatrix}
    1-\alpha & \alpha & 0 \\
    1-\alpha & 0 & \alpha \\
    0 & 0 & 1
\end{pmatrix}$
& $\alpha \in \left[0,\frac{1}{2}\right]$ \\

\bottomrule
\end{tabular}
\caption{Realising matrices for $C_1$--$C_5$} \label{table}
\end{table}

We note that the realising matrices corresponding to the boundary of \(\xi_3\), shown in Table \ref{table}, induce realising matrices for the entire region through the construction in \cite{vagenende2025star}. Figure \ref{fig:possible2} additionally gives the region of all possible real eigenvalue pairs $(\lambda_2,\lambda_3)$ arising from $\mathcal{S}_3$, determined by the constraints $-1 \le \lambda_3 \le \lambda_2 \le 1$ and the JLL-inequality $\mathrm{tr}(S)=1+\lambda_2+\lambda_3 \ge 0$, which is a necessary and sufficient condition  for \(3 \times 3\) stochastic matrices with real spectrum \cite{suleimanova1949stochastic}. We mention that complex eigenvalues can also occur for \(3 \times 3\) stochastic matrices, but that we are only considering real eigenvalues here. In comparison with the stochastic matrices with real spectrum, the region $\xi_3$ represents a substantial contraction.

\section{Reduction Theorem for \(n \geqslant4\)}\label{Reduction theorem}
From Section \ref{Monotone eigenvalue regions}, we already have a complete determination of \(\Xi_1\), \(\Xi_2\) and \(\Xi_3\). For \(n \geqslant 4\) generalised arguments analogous to those in Lemma \ref{properties_dominance_matrix} show that the dominance matrix has column sums less than $1$ and a positive trace. Further, \(\Xi_n\) is star-convex with respect to the origin for \(n \geqslant 4\), so it suffices to determine the boundary. 

Moreover, we prove the Reduction Theorem \ref{ReductionTheorem}. According to this theorem, we can reduce the eigenvalue region \(\Xi_n\) for \(\mathcal{M}_n\) to a subset of the eigenvalue region \(\Theta_{n-1}\) for \(\mathcal{S}_{n-1}\). The proof of this theorem is build on the following lemma.
\begin{lemma}[\cite{minc1988nonnegative}]\label{Minc}
    If \(A\) is an \( n \times n\) non-negative matrix with positive maximal eigenvalue \(r\) and a corresponding positive eigenvector \(x = (x_1, x_2, \ldots, x_n)\), then \( (1/r) \cdot D^{-1}AD \) is a stochastic matrix, with the diagonal matrix \(D = \mathrm{diag}(x_1, x_2, \ldots , x_n)\). Moreover, if \( \sigma(A) = \{\lambda_1, \ldots, \lambda_{n} \}\) and  \( \sigma((1/r) \cdot D^{-1}AD) = \{\mu_1, \ldots, \mu_{n} \}\) are the spectra of respectively \(A\) and \((1/r) \cdot D^{-1}AD\), then
\(\mu_i = \lambda_i/r \text{ for } i \in \{1, \ldots, n\}.\)
\end{lemma}

\begin{theorem}[Reduction Theorem] \label{ReductionTheorem}
    For every \(n \geqslant 4\): \(\Xi_{n} \subseteq \Theta_{n-1} \). 
\end{theorem}
\begin{proof}
Let \(M\) be an \(n \times n\) monotone matrix and \(\sigma(M) = \{1, \lambda_2, \ldots, \lambda_{n} \}\), for \(n \geqslant 4 \), with \(D(M)\) its \((n-1) \times (n-1)\) dominance matrix. It is clear that \(1 \in \Theta_{n-1}\) because \(1\) is an eigenvalue of every \((n-1) \times (n-1) \) stochastic matrix. Our aim is to prove that \(\sigma(M) = \{1, \lambda_2, \ldots, \lambda_{n}\} \subseteq \Theta_{n-1}\), so it suffices to proof that \(\{\lambda_2, \ldots, \lambda_{n} \} \subseteq \Theta_{n-1} \). In order to do so, we use the dominance matrix \(D(M)\) with \(\sigma(D(M)) = \{\lambda_2, \ldots, \lambda_{n} \}\).

The dominance matrix \(D(M)\) can be transformed into a block upper triangular matrix (known as the Perron-Frobenius normal form), so we can assume that \(D(M)\) is a non-negative irreducible matrix. It follows from the Perron-Frobenius theorem for irreducible matrices \cite{seneta1973non} that \(D(M)\) has a strictly positive maximal eigenvalue \(r= \lambda_2\) and a strictly positive maximal eigenvector. It follows from Lemma \ref{Minc} that \(S = (1/\lambda_2) \cdot D^{-1} \cdot D(M) \cdot D \) is an \((n-1) \times (n-1)\) stochastic matrix. If \( \sigma(S) = \{\mu_1, \ldots, \mu_{n-1} \}\), then we have the following link between the eigenvalues of \(S\) and \(D(M)\):
        \[\mu_i = \lambda_{i+1}/\lambda_2 \text{ for } i \in \{1, \ldots, n-1\}.\]

Since \(\lambda_2\) is an eigenvalue of the stochastic matrix \(M\), necessarily holds \(\lambda_2 \leq 1\).

Hence, \(|\mu_i| = |\lambda_{i+1}|/\lambda_2 \geqslant |\lambda_{i+1}| \). Because each \(\mu_i\) lies in the region \(\Theta_{n-1}\) and, moreover, this region is star-convex, it follows that also \(\lambda_{i+1}\) lies in \(\Theta_{n-1}\). Thus \( \sigma(D(M)) = \{\lambda_2, \ldots, \lambda_{n} \} \subseteq \Theta_{n-1} \), which completes the proof.
\end{proof}

\section{Conclusions and further research}\label{Conclusions and further research}

This paper presents a first step in analysing the eigenvalue regions for monotone stochastic matrices. We establish fundamental properties of the dominance matrix $D(M)$ and characterise the conditions under which a non-negative matrix qualifies as a dominance matrix. The eigenvalues of monotone stochastic matrices are examined from two perspectives.

First, viewed individually, we determine the eigenvalue region for $n \times n$ monotone stochastic matrices for $1 \leq n \leq 3$, together with the realising matrices of these regions. Second, viewed collectively, we characterise the set \(\xi_3\) of all possible eigenvalue pairs $(\lambda_{2},\lambda_{3})$ arising from $3 \times 3$ monotone stochastic matrices, again accompanied by corresponding realising matrices. In both approaches, the regions for monotone stochastic matrices are markedly smaller than those obtained for general stochastic matrices. Moreover, we prove a reduction theorem which shows that, for $n \geq 4$, the eigenvalue region for $n \times n$ monotone stochastic matrices is contained within that of $(n-1) \times (n-1)$ stochastic matrices.

The eigenvalue regions for \(\mathcal{M}_n\) are completely determined for \(1 \leq n \leq 3\). However, for \(n \geq 4\), the problem remains unsolved. Although this paper introduces the reduction theorem \ref{ReductionTheorem} to provide additional restrictions, a comprehensive examination is necessary to establish an exact characterisation.

Knowledge about the eigenvalues of monotone stochastic matrices holds potential applications in Markov theory. Consequently, future research should focus on linking the spectral properties of monotone stochastic matrices to their corresponding Markov chains, as eigenvalues offer critical insights into both the long-term behavior and the rate of convergence.

\bibliography{bibl}

@book{Bartholomew,
	Author = {Bartholomew, D.J. and Forbes, A.F. and McClean, S.I.},
	Publisher = {Wiley},
	Series = {Wiley series in probability and mathematical statistics},
	Title = {Statistical techniques for manpower planning. Second edition.},
	Year = {1991}}

@article{seneta1973non,
  title={Non-negative matrices: an introduction to theory and applications},
  author={Seneta, Eugene},
  year={1973}
}

@book{Karpelevich,
	Author = {F.I. Karpelevich},
	Isbn = {978-0-8218-3116-8},
        Volume = {140},
	Pages = {79-100},
	Publisher = {American Mathematical Society Translations},
	Series = {Series 2},
	Title = {Eleven Papers Translated from the Russian: On characteristic roots of matrices with nonnegative elements.},
	Year = {1988}}

@book{Dynkin,
	Author = {Dmitriev, N.A. and Dynkin, E.B.},
	Isbn = {978-0-8218-3116-8},
        Volume = {140},
	Pages = {57-77},
	Publisher = {American Mathematical Society Translations},
	Series = {Series 2},
	Title = {Eleven Papers Translated from the Russian: Characteristic roots of stochastic matrices.},
	Year = {1988}}

@article{guerry2022monotone,
  title={On monotone Markov chains and properties of monotone matrix roots},
  author={Guerry, Marie-Anne},
  journal={Special Matrices},
  volume={11},
  number={1},
  pages={20220172},
  year={2022},
  publisher={De Gruyter}
}

@article{pillai2005perron,
  title={The Perron-Frobenius theorem: some of its applications},
  author={Pillai, S Unnikrishna and Suel, Torsten and Cha, Seunghun},
  journal={IEEE Signal Processing Magazine},
  volume={22},
  number={2},
  pages={62--75},
  year={2005},
  publisher={IEEE}
}

@article{delbianco2023markov,
  title={Markov chains, eigenvalues and the stability of economic growth processes},
  author={Delbianco, Fernando and Fioriti, Andr{\'e}s and Tohm{\'e}, Fernando},
  journal={Empirical Economics},
  volume={64},
  number={3},
  pages={1347--1373},
  year={2023},
  publisher={Springer}
}

@article{racoceanu1995new,
  title={On a new method of Markov chain reduction},
  author={Racoceanu, D and Elmoudni, A and Ferney, M and Zerhouni, S},
  journal={Mathematical Modelling of Systems},
  volume={1},
  number={3},
  pages={199--229},
  year={1995},
  publisher={Taylor \& Francis}
}

@article{meyer2000applied,
  title={Applied Linear Algebra and Matrix Analysis},
  author={Meyer, Carl D and Schmeider, Hans and others},
  journal={Philadelphia, PA, USA: SIAM},
  year={2000}
}

@article{baake2022equal,
  title={On equal-input and monotone Markov matrices},
  author={Baake, Michael and Sumner, Jeremy},
  journal={Advances in Applied Probability},
  volume={54},
  number={2},
  pages={460--492},
  year={2022},
  publisher={Cambridge University Press}
}

@article{jarrow1997markov,
  title={A Markov model for the term structure of credit risk spreads},
  author={Jarrow, Robert A and Lando, David and Turnbull, Stuart M},
  journal={The review of financial studies},
  volume={10},
  number={2},
  pages={481--523},
  year={1997},
  publisher={Oxford University Press}
}

@book{minc1988nonnegative,
  title={Nonnegative matrices},
  author={Minc, Henryk},
  volume={170},
  year={1988},
  publisher={Wiley New York}
}

@article{domka2022spectrum,
  title={On spectrum of Metzler matrices.},
  author={Domka, Micha{\l} and Mitkowski, Wojciech},
  journal={Przeglad Elektrotechniczny},
  volume={98},
  number={12},
  year={2022}
}

@article{mashreghi2007conjecture,
  title={On a conjecture about the eigenvalues of doubly stochastic matrices},
  author={Mashreghi, Javad and Rivard, Roland},
  journal={Linear and Multilinear Algebra},
  volume={55},
  number={5},
  pages={491--498},
  year={2007},
  publisher={Taylor \& Francis}
}

@article{kirkland1992eigenvalue,
  title={An eigenvalue region for Leslie matrices},
  author={Kirkland, Steve},
  journal={SIAM journal on matrix analysis and applications},
  volume={13},
  number={2},
  pages={507--529},
  year={1992},
  publisher={SIAM}
}

@article{conlisk1990monotone,
  title={Monotone mobility matrices},
  author={Conlisk, John},
  journal={Journal of Mathematical Sociology},
  volume={15},
  number={3-4},
  pages={173--191},
  year={1990},
  publisher={Taylor \& Francis}
}

@article{daley1968stochastically,
  title={Stochastically monotone Markov chains},
  author={Daley, Daryl J},
  journal={Zeitschrift f{\"u}r Wahrscheinlichkeitstheorie und verwandte Gebiete},
  volume={10},
  number={4},
  pages={305--317},
  year={1968},
  publisher={Springer}
}

@article{vagenende2025star,
  title={Star-Convexity of the Eigenvalue Regions for Stochastic Matrices and Certain Subclasses},
  author={Vagenende, Brando and Verbeken, Brecht and Guerry, Marie-Anne},
  journal={Mathematics},
  volume={13},
  number={12},
  pages={1--10},
  year={2025},
  publisher={MDPI}
}

@article{loewy1978note,
  title={A note on an inverse problem for nonnegative matrices},
  author={Loewy, Raphael and London, David},
  journal={Linear and Multilinear Algebra},
  volume={6},
  number={1},
  pages={83--90},
  year={1978},
  publisher={Taylor \& Francis}
}

@article{johnson2025perron,
  title={Perron similarities and the nonnegative inverse eigenvalue problem},
  author={Johnson, Charles and Paparella, Pietro},
  journal={Transactions of the American Mathematical Society},
  volume={378},
  number={12},
  pages={8361--8389},
  year={2025}
}

@inproceedings{suleimanova1949stochastic,
  title={Stochastic matrices with real characteristic values},
  author={Suleimanova, HR},
  booktitle={Dokl. Akad. Nauk SSSR},
  volume={66},
  pages={343--345},
  year={1949}
}

\end{document}